\documentclass[10pt]{amsart}
\usepackage[]{amsmath, amsthm, amsfonts}
\usepackage{graphicx}
\usepackage[all]{xy}

\newcommand {\C} {{\mathbb C}}
\newcommand {\R} {{\mathbb R}}
\newcommand {\Z} {{\mathbb Z}}
\newcommand {\Q} {{\mathbb Q}}
\newcommand {\PP} {{\mathbb P}}
\newcommand {\dt} {\bullet}
\newcommand {\F} {{\mathcal F}}
\newcommand {\I} {{\mathcal I}}
\newcommand {\HH} {{\mathbb H}}
\newcommand {\cH} {{\mathcal H}}
\newcommand {\G} {{\mathcal G}}
\newcommand {\LL} {{\mathcal L}}
\newcommand {\A} {{\mathcal A}}
\newcommand {\E} {{\mathcal E}}
\newcommand {\B} {{\mathcal B}}

\newcommand {\PVHS} {\mbox{PVSH}}
\newcommand {\GVHS} {\mbox{GVSH}}

\newtheorem{thm}[subsection]{Theorem}
\newtheorem{cor}[subsection]{Corollary}
\newtheorem{lemma}[subsection]{Lemma}

\newtheorem{remark}[subsection]{Remark}

\begin{document}

\title{The Leray spectral sequence is motivic}

\author{ Donu Arapura}
 \address{
Department of Mathematics\\
 Purdue University\\
 West Lafayette, IN 47907\\
 U.S.A.}
 \thanks{Author partially supported by the NSF}

\email{dvb@math.purdue.edu}

\maketitle

Our goal is  to prove that the Leray spectral sequence
associated to a map of algebraic varieties is motivic in the following sense:
 If the singular cohomology groups of the category of   
quasiprojective  varieties defined over a subfield of $\C$
 can be canonically endowed with some additional
structure, such as its  mixed Hodge structure or its motive in Nori's
sense.  Then the Leray spectral
sequence for a projective map  in this category is automatically 
compatible with this structure. 
Establishing the  compatibility of  Leray with mixed Hodge structures was 
our primary objective, and there  were some important precedents for this.
Zucker  \cite{zucker} established compatibility for a map to  a curve.
This was proved following a delicate analysis of
 variations of Hodge structure on curves and their cohomology.
Saito, in his rather formidable work on
 (mixed) Hodge modules  \cite{saito1, saito2}, addressed the
compatibility issue in general, albeit somewhat indirectly.
 His theory mirrors  the theory of perverse $l$-adic sheaves \cite{bbd}. 
Among other things, Saito obtained an analogue of the decomposition theorem,
which gives, in practice, a replacement for the Leray spectral sequence.
In view of this history, it came as a pleasant surprise to author
that one could get quite far by using nothing beyond the standard
functoriality properties established by Deligne \cite{deligne2}.

A very brief outline of the proof of the main result (theorem
\ref{thm:gen}) is as follows. We first show that the spectral sequence 
associated to a filtration by subvarieties is motivic. Then 
we identify the  Leray spectral sequence 
with  such a spectral sequence associated to a suitable filtration. 
The quasiprojectivity assumption
is imposed mainly for convenience since it permits considerable simplifications
of the arguments.
When the map is  smooth and the coefficients are rational, the
Leray spectral sequence degenerates, and we can exhibit  its motivic  
structure in a more elementary fashion by theorem \ref{thm:smooth}. 
This special case is already sufficient
for many purposes.

As an antidote to all these generalities, we conclude the paper with a
simple geometric example to illustrate how these results can be
applied (some readers may wish to start with this). 
These results can be used to simplify or justify arguments that have
already appeared in the literature  (e.g. \cite[6.3.1]{as},
\cite[3.4]{lewis}).
Another byproduct is a relatively cheap construction of  a mixed Hodge
structure on the cohomology of variation of Hodge structure of
geometric origin. This can be used as an alternative to
Saito's machinery for certain applications (e.g. \cite{hain}).
We have not checked compatibility with Saito's construction, 
but we expect that it should be true.

I would like to express my warmest thanks to Madhav Nori for sharing
his insight into these problems with me. The special case
(\ref{thm:smooth}) was found by him independently, and  the proof of
the general case (\ref{thm:gen}) was based on his observation that  it
is possible to construct  ``cellular decompositions'' for affine
varieties with respect to  constructible sheaves.  I would also like to 
thank Chris Peters for pointing out a gap in the original proof
of theorem~\ref{thm:smooth}, and for suggesting  a remedy; I will say more
about that below.
My thanks also to
the referee for suggesting some improvements in the exposition.

\section{Preliminaries}

Let us fix some conventions and terminology.
We will drop indices and other
decorations from the notation whenever  they seem 
unnecessary. 
Let $V_\kappa$ be the category of quasiprojective 
schemes over  be a field $\kappa$. We refer to these objects as
varieties whether not they are irreducible. The dimension
is the dimension of the largest component.
Let $V^2_\kappa$ be the
category whose objects are pairs $(X,Y)$ with
$Y\subseteq X$ closed. A morphism from $(X,Y)\to (X',Y')$
is a morphism of $\kappa$-schemes $X\to X'$ such that $f(Y)\subseteq
Y'$. The category $V^3_\kappa$ of triples of
 closed embeddings $Z\subseteq Y\subseteq X$ in $V_\kappa$ can be defined in a similar
fashion. We define a {\em weak cohomology} theory over $\kappa$ to be a
 functor $\cH:V_\kappa^2\to \A$ to an Abelian category such that
for any object $Z\subseteq Y\subseteq X$ in $V^3_\kappa$,
there are natural transformations
$$\cH^i(Y,Z)\to \cH^{i+1}(X,Y)$$
which fit into a long exact sequence:
\begin{equation}\label{eq:longexact}
\ldots \cH^i(Y,Z)\to \cH^{i+1}(X,Y)\to \cH^{i+1}(X,Z)
\ldots 
\end{equation}
As usual, we write $\cH^i(X) = \cH^i(X,\emptyset)$.
$\A$ will be called the coefficient category.

Here is a  list of examples.

\begin{enumerate}

\item[(Betti)]
  Let $\kappa=\C$, then take $\cH=H$ to be singular (or Betti) cohomology 
with $\Z$ coefficients.  In this example
 the coefficient category is the category of Abelian groups $Ab$

\item[(MHS)]  Let $\kappa=\C$, take $\cH = H_{MHS}$ as $H$ together with its
canonical mixed Hodge structure \cite{deligne2}.
The coefficient category  is the category
of mixed Hodge structures.

\item[($\sigma$-Betti)]
Let $\sigma:\kappa\hookrightarrow \C$ be an embedding.
The composition $H_\sigma$ of $H$ with base extension  by $\sigma$
gives a weak theory with coefficients $Ab$.
one can define $H_{\sigma,MHS}$ in a similar fashion.

\item[(l-adic)] Let $l$ be a prime different from the characteristic
of $\kappa$, let $\bar \kappa$ denote the separable closure and let
$\cH_l^i(X,Y) = H_{et}^i(X\otimes \bar \kappa, j_{Y,X!}\Q_l)$
where $j_{X,Y}: X-Y\hookrightarrow X$ is the inclusion.
The long exact sequence (\ref{eq:longexact}) is induced from
$$0\to j_{X,Y!}\Q_l\to  j_{X,Z!}\Q_l\to
 j_{Y,Z!}\Q_l \to 0 $$
$\A$ is the category of $l$-adic representations of 
$Gal( \bar \kappa/\kappa)$.

\item[(AH)] When $\kappa$ is embeddable into $\C$,
Jannsen \cite[6.11.1]{jannsen} has essentially constructed
a weak cohomology theory $H_{AH}$ with values in his category of 
integral mixed realizations $MR_k$. Here we ignore Tate twists. 
The cohomology of a pair 
can defined by combining Jannsen's methods with the simplicial mapping
cone construction \cite[6.3]{deligne2}.

\item[(ECM)] When $\kappa$ is embeddable into $\C$, Nori (unpublished) has constructed
the universal weak cohomology theory $H_{ECM}$. The coefficient
category is  his category of effective cohomological motives. 

\end{enumerate}

We assume from now on that $\kappa$ is embeddable into $\C$.
Fix one such  embedding  $\sigma:\kappa\hookrightarrow \C$. We will call
a weak cohomology theory $\cH:V_\kappa^2\to \A$ an {\em enriched
 (rational) $\sigma$-Betti} theory if there exist a faithful exact functor
$\Phi:\A\to Ab$ ( $\Phi:\A\to \Q\mbox{-(vector spaces)}$)
such that $ H_\sigma = \Phi\circ \cH$ ( $H_\sigma\otimes\Q =\Phi\circ \cH)$).
Both $H_{\sigma,MHS}$ with the obvious $\Phi$, and $H_{AH}$
and $H_{ECM}$ provided examples of enriched Betti theories.

By a spectral sequence in an Abelian
category $\A$,  we mean a collection of bigraded objects $E_r^{\dt\dt}$
with differentials of bidegree $(r,-r+1)$, for $r\ge r_0$ for some
$r_0\ge 0$,
together with a filtered graded object $(H^\dt, F^\dt)$ with
an isomorphisms
 $$E_{r+1}\cong H^\dt(E_r),\  E_\infty \cong Gr_F H.$$
We will denote the whole thing by $E_{r_0}\Rightarrow H$ or simply
$E_{r_0}$. Clearly the set of spectral
sequences (for fixed $r_0$) form a category $SS(\A,r_0)$ in its own right.
Given an exact functor of Abelian categories $\Phi:\A\to \B$ and a
spectral sequence $E_{r_0}$ in $\B$, by a lift to $\A$ we mean
a spectral sequence $\E_{r_0}$ in $\A$ and an isomorphism $E_{r_0}\cong 
\Phi(\E_{r_0})$. More generally, given a functor $F:C\to SS(\B,r_0)$,
a natural lift is a functor $F':C\to SS(\A,r_0)$ and
a natural isomorphism $\Phi\circ F' \cong F$.

Given a continuous map of topological spaces $f:X\to Y$ and sheaf
$\F$ on $X$, the Leray spectral sequence 
\begin{equation}\label{eq:LerayDef}
\LL_2^{p,q}(f,\F) = H^p(Y,R^q f_*\F)
\end{equation}
$$\LL_\infty^{p,q}=Gr_L^pH^{p+q}(X,\F)$$
provides the basic example of a spectral sequence in $Ab$. 
The  Leray filtration is given by
$$
L^pH^i(X,\F) = image[H^i(Y,\tau_{\le i-p} \R f_*\F)\to H^i(Y,\R 
f_*\F)
\cong H^i(X,\F)].$$
This is functorial in  the sense that a commutative diagram
$$
\xymatrix{
 X'\ar[r]^{g}\ar[d]^{f'} & X\ar[d]^{f} \\
 Y'\ar[r]^{h} & Y
}
$$
induces morphism 
\begin{equation}\label{eq:basechangeLeray}
\LL_2^{pq}(f,\F)\to \LL_2^{pq}(f', g^*\F).
\end{equation}

Functoriallity follows from the construction which  we briefly recall
along with  some related homological algebra.
A biregular filtered complex $(K^\dt, F^\dt)$ in
an Abelian category gives rise to
a spectral sequence
\begin{equation}\label{eq:specseq}
E_1^{pq}(K,F)=H^{p+q}(Gr_F^pK^\dt)\Rightarrow H^{p+q}(K^\dt)
\end{equation}
This can be constructed using explicit formulas or
in terms of exact couples which we recall later.
 The shifted filtration (``filtration decal\'e'') \cite[1.3.3]{deligne2} is
\begin{equation}\label{eq:Dec1}
Dec\, (F)^pK^n = \{\alpha\in F^{p+n}K^n\,|\, d\alpha\in F^{p+n+1}\}.
\end{equation}
For example, the canonical filtration $\tau^\dt= \tau_{\le -\dt}$ is just
$Dec\, (G)$, where $G$ is the trivial filtration $G^0=K$, $G^1 = 0$.
The spectral sequence with respect to the shifted and unshifted
filtrations are the same up to a change of indices \cite[1.3.4]{deligne2}:
\begin{equation}\label{eq:Dec}
 E_r^{pq}(K, Dec\, (F))\cong E_{r+1}^{2p+q, -p}(K, F)
\end{equation}

A morphism of filtered complexes is a filtered
quasiisomorphism if it induces an isomorphism on the associated graded
complexes.  Such a morphism induces an isomorphism of (\ref{eq:specseq}).
The  filtered derived category \cite[7.1, 1.4]{deligne2} is obtained from the
category of filtered complexes by inverting filtered
quasiisomorphisms; the spectral sequence (\ref{eq:specseq}) extends
to this setting. Any object $\F$ of the derived category, gives
rise to a well defined object $(\F,\tau)$ of the filtered derived category.
A filtered acyclic resolution is a filtered 
quasiisomorphism to a filtered complex  whose associated graded
is acyclic. Derived functors, which exist in the filtered context,
 are computed by applying the functor to a filtered acyclic
 resolution of the original complex. 
In particular, given a filtered complex of sheaves
$(K^\dt, F^\dt)$, we have  an object $\R\Gamma(K,F)$ of the filtered
derived category of Abelian groups, and this maps to $\R\Gamma(K)$
under the forgetful functor to the ordinary derived category.
 There is a spectral sequence
\begin{equation}\label{eq:specseq2}
E_1^{pq}(\R\Gamma(K,F))= \HH^{p+q}(Gr_F^pK^\dt) \Rightarrow \HH^{p+q}(K^\dt)
\end{equation}
Suppose $f:X\to Y$ is a continuous map, and $\F$ a sheaf
on $X$. The canonical filtration 
 on $\R f_*\F$ induces the Leray spectral
sequence with a shift of indices \cite[1.4.8]{deligne2}:
\begin{equation}\label{eq:shiftEpqr}
  E_r^{pq}(\R\Gamma(\R f_*\F, \tau))\cong \LL_{r+1}^{2p+q, -p}(f,\F)
\end{equation}  

\section{Smooth Projective Maps}

For the remainder of this paper, we will be concerned with varieties
defined over a subfield $\sigma:\kappa\hookrightarrow \C$.
By a sheaf on a such variety $X$, we will mean a sheaf on
$(X\times_\sigma\, Spec\,\C)^{an}$. We write $H^i(X,\F)$
instead of $H^i((X\times_\sigma\, Spec\,\C)^{an},\F)$.
We say that a  sheaf $\F$ of $\Z$-modules  on  $X$
is weakly constructible  if there is 
a partition of $X$ into Zariski locally closed sets $\mathcal{Z}$
defined over $\kappa$  such that
$\F|_{Z}$ is locally constant for each $Z\in \mathcal{Z}$; $\F$
is constructible if in addition the stalks are finitely generated.

\begin{thm}\label{thm:smooth}
Let $Y$ be a smooth complex  irreducible  quasiprojective variety.
Then  for each $p$, there  a morphism
$Y_p\to Y$ with the following property:
If $f:X\to Y$ is   a smooth projective morphism,
then
$$L^pH^i(X,\Q) = ker[H^i(X,\Q)\to H^i(X_p,\Q)]$$
where $X_p = f^{-1}Y_p$.
There are morphisms $Y_{p-1}\to Y_p$ over $Y$ compatible with
inclusions $L^p\subseteq L^{p-1}$.
If $ Y$ is defined over an algebraically closed subfield 
$\kappa\subseteq \C$, then $Y_p$ etcetera can be  defined
over $\kappa$.
\end{thm}

\begin{remark}
 In an earlier incarnation of this paper, we stated and attempted to
prove this theorem without assuming smoothness of $Y$. However,
C. Peters pointed a gap in the original form of
lemma~\ref{lemma:lef}, and suggested a fix using work
of \cite{hamm} among things; $H$ would have to be
 allowed to be replaced by a multiple. In the interests of simplicity,
we opted to state and prove this under a more restrictive hypotheses. 
In any case, this theorem is superceded by theorem~\ref{thm:gen}.
\end{remark}

\begin{lemma}[Jouanolou]\label{lemma:jou}
Let $Y$ be a quasiprojective variety over a field $\kappa$. There exists 
an affine $\kappa$-variety $Y'$ and a morphism $Y'\to Y$ where the
fibers are isomorphic to an affine space.
\end{lemma}

\begin{proof}
See \cite[1.5]{jou}
\end{proof}

Call a $\kappa$-variety $F$ contractible, if $(F\times \C)^{an}$ is 
contractible.

\begin{cor}\label{cor:jou}
  Let $X\to Y$ be a morphism of quasiprojective varieties over
  $\kappa$. There exists a commutative diagram
$$
\xymatrix{
 X'\ar[r]\ar[d]& X\ar[d] \\
 Y'\ar[r] & Y
}
$$
where both $X'$ and $Y'$ are affine, and the fibers of the horizontal
maps are contractible.
\end{cor}

\begin{proof}
  Apply the lemma to $Y$ to obtain $Y'$, and then again
to $X\times_{Y}Y'$ to obtain $X'$.
\end{proof}

We need the following version of the weak Lefschetz theorem.

\begin{lemma}\label{lemma:lef}
 Let $L$ be a locally constant sheaf on  an $n$
  dimensional  nonsingular complex affine variety $Y$. Then
$H^i(Y, L) = 0$ for $i>n$,
 and $H^i(Y,L)\to H^i(H, L)$ is injective for $i<n$ and any
general affine hyperplane section $H$.
\end{lemma}

\begin{proof}
The vanishing follows from the topological analogue of Artin's theorem 
about the vanishing of the cohomology of constructible
sheaves on $Y$ in degrees $>n$. (A proof can be found in \cite{gm}.
One can also deduce this from lemma~\ref{lemma:beil} and induction.)
For the second statement, we apply \cite[thm 6.1.1]{as}
\end{proof}

\begin{proof}[Proof of theorem \ref{thm:smooth}]
We first treat the special case where $Y$ is affine, and prove this by induction
on dimension. The  initial case where $Y$ is a point, or more
generally a union of points, is vacuously true. We can take $Y_0
=\emptyset$.  
Assume that $\dim Y = n>0$. Choose a general hyperplane section $H$
(defined over $\kappa$), and consider the cartesian square
$$
\left.
\begin{array}{ccc}
 X_H & \to  & X \\
 g\downarrow &  & \downarrow f\\
 H & \to & Y
\end{array}
\right.
$$
This induces a morphism of Leray spectral sequences
$$
 \LL_2^{pq}(f,\Q)\to \LL_2^{pq}(g,\Q)
$$
Furthermore, by \cite{deligne}, both spectral sequences degenerate
at $E_2$. Together with lemma \ref{lemma:lef}, this implies that 
$\LL_\infty^{pq}(f,\Q)$ and $ \LL_\infty^{pq}(g,\Q)$ vanish for $p>n$ and
$p>n-1$ respectively. Moreover, the map
$$\LL_\infty^{pq}(f,\Q)\to \LL_\infty^{pq}(g,\Q)$$
is an injective if $p<n$.
 Thus the smallest  step of the Leray filtration
$$L^nH^i(X,\Q) = \LL_\infty^{n,i-n}(f,\Q) =
 ker[H^i(X, \Q) \to H^i(X_H,\Q)]$$
Therefore, we can take $Y_n = H$. Suppose that $p< n$. By induction,
there is a morphism of varieties $H_p\to H$ such that
$$L^pH^i(X_H,\Q) = ker[H^i(X_H,\Q)\to H^i(f^{-1}H_p,\Q)]$$
Since $L^p$ is compatible with the restriction $H^i(X,\Q)\to
H^i(X_H,\Q)$,
it  follows that
$$L^pH^i(X,\Q)\subseteq  ker[H^i(X,\Q)\to H^i(f^{-1}H_p,\Q)]$$
We claim equality holds. Suppose to the contrary, that there 
exists  $\alpha$ in the complement of the right and left hand
sides. Choose the least $q$ such that
$\alpha\in L^qH^i(X,\Q)$. We have $q<p$.
Then the class $Gr_L^q(\alpha)\not= 0$, which implies
that $Gr_L^q(\alpha|_{X_H})\not= 0$. Therefore $\alpha|_{X_H}\notin
L^{p}$. It follows that its restriction to $f^{-1}H_p$ is nonzero,
and this is a contradiction.  This proves the claim. Setting $Y_p = H_p$
completes the proof for affine $Y$.

Let $Y'\to Y$ be as in lemma \ref{lemma:jou} and $X' = X\times_Y Y'$.
Since the fibers of $Y'\to Y$ are contractible, the maps
$f^{an}$ and ${f'}^{an}:{X'}^{an}\to {Y'}^{an}$ are homotopy equivalent. Thus
there is an isomorphism $H^*(X)\cong H^*(X')$ strictly compatible with 
Leray filtrations, i.e. that the induced map on $Gr_L$
is an isomorphism. By the previous case, there exists
$Y_p'\to Y'$ such that $L^pH^i(X')$ is the kernel of the
map to $H^i({f'}^{-1}Y_p')$. Then the composition $Y_p=Y_p'\to Y$
will have the required property.
\end{proof}

\begin{cor}\label{cor:MHSleray}
The Leray filtration on $H^i(X,\Q)$ is a filtration by sub-mixed
Hodge structures.
More generally, for any enriched  rational $\sigma$-Betti theory
$(\cH:V_\kappa^2\to \A,\Phi)$, the Leray filtration on $H^i(X,\Q)$
is the image of a filtration on $\cH^i(X)$.
\end{cor}

\section{Main theorem}

Fix an embedding  $\sigma:\kappa\hookrightarrow \C$
and   an enriched $\sigma$-Betti theory
$(\cH:V_\kappa^2\to \A,\Phi)$. 
Let $f:X\to Y$ be a projective morphism
of quasiprojective varieties over $\kappa$,
then we can form the  Leray spectral sequence
\begin{equation}\label{eq:Leray}
E_2 = H_\sigma^p(Y,R^qf_{\sigma*}\Z) \Rightarrow
H_\sigma^{p+q}(X,\Z)
\end{equation}
where $f_\sigma$ is the analytic map associated to 
$f\times_\sigma\, Spec\,\C$.
As explained in section 1, a natural lift of (\ref{eq:Leray}) consists
of a spectral sequence
$$\E_2^{pq}(f)\Rightarrow \cH$$
in $\A$ which is  compatible with base change (\ref{eq:basechangeLeray})
 and maps to (\ref{eq:Leray}) under $\Phi$.

\begin{thm}\label{thm:gen}
The Leray spectral sequence  has a  natural lift to $\A$.
(Thus it is motivic.)
\end{thm}

\begin{cor}\label{cor:LerayMHS}
The Leray spectral sequence of a projective map
of complex quasiprojective varieties can be lifted naturally to the
category of mixed Hodge structures.
\end{cor}

Since the proof of the theorem will be spread out over a long  series  of 
lemmas, it may be worthwhile  to outline the main steps. (We suppress $\sigma$
in the outline.)

\begin{enumerate}
\item[A)] If $X$ has an increasing filtration $X_\dt$ by closed sets, then
 lemma~\ref{lemma:mhs} will yield a spectral sequence of Abelian groups
$$E_1^{pq}(X_\dt, \Z)\Rightarrow H^{p+q}(X,\Z)$$
which is motivic in the sense that it has a natural lift to $\A$.

\item[B)] If $X_\dt$ is given by pulling back a 
  filtration $Y_\dt$ from $Y$, then in lemma~\ref{lemma:leray}
 we will construct a morphism
$$ H^p(Y,R^qf_*\Z)\to E_2^{pq}(X_\dt,\Z)$$
from the Leray spectral sequence to the one above
from $E_2$ onwards. Note that this morphism is given in the
category of Abelian groups; the left side has no motivic
structure a priori.

\item[C)] As in the proof of theorem~\ref{thm:smooth}, we can reduce
to the case where $Y$ is affine. Then by 
applying lemma~\ref{lemma:existenceofcell} to the direct images $R^qf_*\Z$,
we can choose the $Y_\dt$ so that the morphism of step B is an isomorphism.
This gives the desired motivic structure on Leray.

\item [D)] The final step is to check independence of choices and
naturality.
\end{enumerate}

It is worth remarking that there certain parallels between this proof 
and the proof of theorem~\ref{thm:smooth}. For example, 
the collection  of varieties $\{Y_\dt\}$ 
play essentially the same role in both proofs.

In order to flesh out the above sketch, we need some more notation.
We can form a category $FV_\kappa$ as follows. The objects
are pairs consisting of a variety $X$  in $V_\kappa$ and an increasing
exhaustive  filtration  
$$ X= \ldots X_{M+1}= X_M\supseteq X_{M-1}\supseteq\ldots X_{-1} = \emptyset$$
by closed sets. A morphism of $f:(X,X_\dt)\to (Y,Y_\dt)$ is a 
morphism of varieties such that $f(X_i) \subseteq Y_i$. 
Suppose that $(Y,Y_\dt)$ is an object of $FV_\kappa$.
 Let $Y_a^o = Y_a-Y_{a-1}$, and let
$j_a:Y_a^o\hookrightarrow Y_a$, $k_a:Y-Y_a\to Y$ and 
$i_a:Y_a\to Y$ denote the inclusions.
Fix a sheaf $\F$ on $Y$, and let $\F_a = j_{a!}\F|_{Y_a^o}$,
$S^a(Y_\dt,\F) = k_{(a-1)!}k_{a-1}^*\F$, $S^0(Y_\dt,\F) = \F$.
We will suppress one or both arguments of $S$ if these are clear from
 context.
We have a filtration by ``skeleta'': 
$\F=S^0(\F)\supset S^1(\F)\supset\ldots$
such that $S^a/S^{a+1}\cong i_{a*}\F_a$. This isomorphism
follows from the snake lemma applied to the diagram
$$
\left.
\begin{array}{ccccccccc}
  &  &  &  &  &  & 0 &  &  \\
  &  &  &  &  &  & \downarrow &  &  \\
  &  & 0 &  &  &  & i_{a*}\F_a &  &  \\
  &  & \downarrow &  &  &  & \downarrow &  &  \\
 0 & \to & S^{a+1} & \to & \F & \to & \F|_{Y_a} & \to  & 0 \\
  &  & \downarrow &  & || &  & \downarrow &  &  \\
 0 & \to & S^{a} & \to & \F & \to & \F|_{Y_{a-1}} & \to & 0 
  
\end{array}
\right.
$$

Given a complex of sheaves of Abelian groups $\F^\dt$,
$S^\dt(Y_\dt, \F^\dt)$ is a filtration by subcomplexes.
The spectral sequence (\ref{eq:specseq2}) for $(\F^\dt, S(Y_\dt,\F^\dt))$
becomes 
\begin{equation}\label{eq:specseqYdt}
E_1^{ab}(Y_\dt,\F^\dt)=  \HH^{a+b}(Y,i_{a*}\F_a^\dt)
\Rightarrow \HH^{a+b}(Y,\F^\dt) 
\end{equation}
The differential
$$d_1:E_1^{a,b}(Y_\dt,\F^\dt)\to E_1^{a+1,b}(Y_\dt,\F^\dt)$$
is given explicitly as a composition of
$$\HH^{a+b}(Y,\F_{a}^\dt)\to \HH^{a+b}(Y,\F^\dt|_{Y_a})$$
and the connecting map
$$\HH^{a+b}(Y,\F^\dt|_{Y_a})\to \HH^{a+b+1}(Y,\F^\dt|_{Y_{a+1}})$$

Given a morphism $f$ and a complex $\F^\dt$ on $Y$, we get
 a morphism of spectral sequences
$$f^*:E_1^{ab}(Y_\dt,\F^\dt) \to E_1^{ab}(X_\dt,f^*\F^\dt).$$
In particular, 
$$(Y,Y_\dt)\mapsto E_1^{ab}(Y_\dt,\Z)$$
is a functor from $FV_\kappa\to SS(Ab,1)$
                         
Since the functors $Gr_{S^a}=(i_a\circ j_a)_!j_a^*$ are exact, a
quasiisomorphism  $\F^\dt\to \G^\dt$  gives rise to a filtered quasiisomorphism
$(\F^\dt,S^\dt(\F^\dt))\to (\G^\dt,S^\dt(\G^\dt))$.
Thus an object $\F$ in the derived category $D^+(Y)$ 
gives rise to well defined
object $(\F, S^\dt(\F))$ in the filtered derived category $DF^+(Y)$.
If $\F$ is replaced by an injective resolution $\I^\dt$, then this is 
represented by $(\I^\dt, S^\dt(\I^\dt))$ thanks to:

\begin{lemma}\label{lemma:Sinjective}
    If $\I$ is injective then $Gr_{S^a}\I$ and $S^\dt(Y_\dt,\I)$ are flasque.
\end{lemma}   

\begin{proof}
    Both $j^*$ and $j_{!}$  preserve injectivity for any open
immersion $j$. Injective sheaves are flasque, and $i_*$ takes flasque
sheaves to flasque sheaves
for any  closed immersion $i$. This proves that $Gr_{S^a}\I$ is
flasque. Since the class of flasque sheaves is closed under
extensions, $S^\dt(Y_\dt,\I)$ is also flasque.
\end{proof}    

We say that an object $(Y,Y_\dt)$ is {\em cellular} with respect
to a sheaf $\F$ if  $H^i(Y,\F_a) = 0,$ for $i\not=a$.
The terminology is suggested by the analogy with $CW$ or cell
complexes in topology with its filtration by skeleta.

\begin{lemma}\label{lemma:cell}
 Suppose that   $H^i(Y,\F_a) = 0$ for $i\not=a$, then
 $$H^i(Y,\F) \cong H^i(E_1^{\dt0}(Y_\dt,\F), d_1)$$
\end{lemma}

\begin{proof}
This follows immediately from the spectral sequence (\ref{eq:specseqYdt}).
\end{proof}

\begin{lemma}[Beilinson]\label{lemma:beil}
If $\F$ is a weakly constructible sheaf on  an $n$
dimensional affine variety $Y$,
there exists a nonempty open set $j:U\hookrightarrow X$ such that 
$$H^i(Y, j_!j^*\F) = 0$$
unless $i = n$.
\end{lemma}

\begin{proof} See \cite{nori}.
\end{proof}

\begin{cor}\label{cor:beil}
Suppose that $\F_1,\F_2,\ldots$ is a finite collection of weakly constructible
sheaves on $Y$ and that $Y' \subset Y$ is a proper closed set. There
exists a   nonempty open set $j:U\hookrightarrow X-Y'$
such that for all $k$,
$$H^i(Y, j_!j^*\F_k) = 0$$
unless $i = n$.
\end{cor}

\begin{proof}
Apply the lemma to $(V,\F_1\oplus \F_2\oplus\ldots)$, where   $V$ is
 an affine subset of $X-Y'$,  
\end{proof}

The existence of cellular filtrations is given by:

\begin{lemma}\label{lemma:existenceofcell}
Let $Y$ be an $n$ dimensional  affine variety.
Suppose that $\F,\ldots$ is a finite collection of weakly constructible
sheaves on $Y$ and that 
$$Y_0'\subset Y_1'\subset\ldots Y_n'= Y$$
 is a chain of closed sets  such  $\dim Y_i'=i$.
There exists a filtration 
$$Y_0\subset Y_1\subset\ldots Y_n= Y$$
 such that
\begin{enumerate}
\item  $\dim Y_i = i$.
\item $Y_i'\subseteq Y_i$.
\item $Y_\dt$ is cellular with respect to $\F,\ldots$
\end{enumerate}
\end{lemma}

\begin{proof}
  This follows from corollary \ref{cor:beil} and induction on $\dim Y$.
\end{proof}

\begin{lemma}\label{lemma:mhs}
    For any pair $(X,X_\dt) $ in $FV_\kappa$,
$$E_1^{pq}(X_\dt,\Z)\Rightarrow H^{p+q}(X,\Z)$$
can be lifted naturally to $\A$.
\end{lemma}

\begin{proof}
The spectral sequence $E_1(X_\dt,\Z)$  can be described using
the language of exact couples \cite[5.9]{weibel}.
We set 
$$D_1 = \bigoplus H^\dt(X, S^\dt(\Z)) =  \bigoplus H^\dt(X,X_{\dt};\Z)$$ 
and 
$$E_1=\bigoplus H^\dt(X,  S^\dt(\Z)/S^{\dt+1}(\Z)) = \bigoplus H^\dt(X_\dt,X_{\dt+1};\Z)$$
with appropriate bigradings. Then we have an exact couple
of Abelian groups

$$
\xymatrix{
 D_1\ar[rr]^{\alpha} &  & D_1\ar[ld]^{\beta} \\
  & E_1\ar[lu]^{\gamma} & 
}
$$
where the morphisms are  obtained from the long exact
sequence associated to the triple $(X,X_\dt,X_{\dt+1})$.
We can define an exact couple in $\A$ with
$${\mathcal D}_1 =   \bigoplus \cH^\dt(X,X_{\dt})$$
$$\E_1= \bigoplus \cH^\dt(X_\dt,X_{\dt+1})$$
in the same way as above. This  exact couple maps to the previous one 
under $\Phi$.
The couple yields  a spectral sequence $\E_1$ in $\A$
which maps to $E_1$ under $\Phi$, and converges to
$\cH^\dt(X)$ with the filtration
$$image[\cH^i(X,X_\dt)\to \cH^i(X)]$$
The exact couple, and hence the spectral sequence, is clearly
functorial in $(X,X_\dt)$.
\end{proof}
   
Suppose that $f:X\to Y$ is a morphism of varieties, and that
$Y_\dt$ is a filtration as above. 
Let $X_\dt = f^{-1}Y_\dt$, then $f:(X,X_\dt)\to (Y,Y_\dt)$ becomes
a morphism.

\begin{lemma}\label{lemma:functorialityofS}
  With the above notation and assumptions,
$\R f_*(\F, S(X_\dt)) \cong (\R f_* \F, S(Y_\dt))$.
for any sheaf $\F$ on $X$. 
\end{lemma}

\begin{proof}
Let
$X_a^o= f^{-1}Y_a^o$, and let
$J_a:X_a^o\hookrightarrow X_a$ and $I_a:X_a\to X$
$K_a:X-X_a\hookrightarrow X$  denote the
inclusions.
 The lemma follows from the  natural equality of sheaves
$$f_*S^{a}(X_\dt, \I) =
k_{a-1!}f_*K_{a-1}^* \I = S^{a}(Y_\dt, f_*\I)$$
for $\I$ on $X$.
\end{proof}

\begin{cor}\label{cor:func1}
    There is an isomorphism of spectral sequences
$$E_1^{pq}(Y_\dt, \R f_*\F)\cong E_1^{pq}(X_\dt, \F)$$
\end{cor}

\begin{cor}\label{cor:func2}
    Suppose that $Y_\dt$ is cellular with respect to 
    $R^bf_*\F$ for all $b$.
    Then there is an isomorphism of bigraded Abelian groups
    $$E_2^{a0}(Y_\dt,R^bf_*\F) \cong  E_2^{ab}(X_\dt, \F)$$
\end{cor}  

\begin{remark}
This will not be an isomorphism of {\bf differential} bigraded groups.
\end{remark}

\begin{proof}
    By the above assumptions, the maps
    $$R^bf_*\F[-b]\leftarrow \tau_{\le b}\R f_*\F \to \R f_*\F$$
    induce isomorphisms
    $$E_1^{a0}(Y_\dt,R^bf_*\F) \cong E_1^{pq}(Y_\dt, \R f_*\F)$$
which will compatible with the differentials.
\end{proof}    

\begin{lemma}\label{lemma:leray}
Let $(Y,Y_\dt)$ be an object of $FV_\kappa$, 
$f:X\to Y$ be a morphism of varieties, and $X_\dt = f^{-1}Y_\dt$.
 For each  sheaf $\F$ on $X$, let
 $$\LL_2^{pq}(f,\F)= H^p(Y,R^qf_*\F)$$
 denote the Leray spectral sequence.
 Then there is a natural  morphism of spectral sequences 
 $$\LL_2^{pq}(f,\F) \to E_2^{pq}(X_\dt, \F)$$
Naturallity means that given objects and morphisms
$(Y',Y_\dt')$, $f':X'\to Y'$,  $X'_\dt = {f'}^{-1}Y_\dt'$,  
and a commutative diagram
$$
\xymatrix{
 (X',X_\dt'),\ar[r]^{g}\ar[d]^{f'} & (X,X_\dt)\ar[d]^{f} \\
 (Y,Y'_\dt)\ar[r]^{h} & (Y,Y_\dt)
}
$$
the induced diagram
$$
\xymatrix{
 \LL_2^{pq}(f', g^*\F)\ar[r]\ar[d] &  \LL_2^{pq}(f,\F)\ar[d] \\
 E_2^{pq}(X'_\dt, g^*\F)\ar[r] &  E_2^{pq}(X_\dt, \F)
}
$$
commutes.
\end{lemma}

\begin{proof}
    Let $\I^\dt$ be an injective resolution of $\F$.
 Set $\Sigma^a =  S^a(X_\dt,\I^\dt)$,    then 
 $(\I^\dt, \Sigma)$ is a filtered acyclic resolution of
    $(\F,S^\dt(X_\dt,\F) )$ by lemma \ref{lemma:Sinjective}.
 If $G$ denotes the trivial filtration $G^0=f_*\I$ and $G^1=0$, then
 the inclusion $(f_*\I,G)\to (f_*\I, f_*\Sigma)$ is compatible with
 the filtrations.   Therefore, we get a map of filtered complexes
 $(f_*\I,\tau)\to (f_*\I, Dec\, (f_*\Sigma))$, see  (\ref{eq:Dec1}).
By definition, there is an exact sequence
$$0\to Dec\, (\Sigma)^p\I^n\to \Sigma^{p+n}\I^n\to 
\Sigma^{p+n}\I^{n+1}/\Sigma^{p+n+1}\I^{n+1}$$
which yields the  exact sequence
\begin{equation}\label{eq:DecSigma1}
0\to f_*Dec\, (\Sigma)^p\I^n\to f_*\Sigma^{p+n}\I^n\to 
f_*(\Sigma^{p+n}\I^{n+1}/\Sigma^{p+n+1}\I^{n+1})
\end{equation}
Furthermore,
\begin{equation}\label{eq:DecSigma2}
f_*(\Sigma^{p+n}\I^{n+1}/\Sigma^{p+n+1}\I^{n+1}) = 
f_*(\Sigma^{p+n}\I^{n+1})/f_*(\Sigma^{p+n+1}\I^{n+1})
\end{equation}
since $(\I^\dt, \Sigma)$ is a filtered acyclic resolution.
Sequences (\ref{eq:DecSigma1}) and (\ref{eq:DecSigma2}) imply that
 $Dec\, ^\dt(f_*\Sigma) = f_*(Dec\, ^\dt \Sigma)$.
 Therefore, 
 we get a map 
 $(f_*\I,\tau)\to (f_*\I, f_*Dec\, (\Sigma))$,
which induces a morphism of spectral sequences
 $$E_1^{pq}(\R\Gamma(f_*\I,\tau))\to
 E_1^{pq}(\R\Gamma(f_*\I, f_*(Dec\,  \Sigma)))$$
The spectral sequence on the left is $\LL$ up to a shift by
(\ref{eq:shiftEpqr}). While the one on the right coincides
with $E_2^{pq}(X_\dt, \F)$ by the same shift thanks to
(\ref{eq:Dec}).

Naturallity follows from lemma \ref{lemma:functorialityofS} and
the commutativity of the diagram
$$
\xymatrix{
 (f_*\I,G)\ar[r]\ar[d] & (f_*\I, S^\dt(Y_\dt, f_*\I))\ar[d]\ar[rd] &  \\
 (f_*g_*g^*\I,G)\ar[r] & (f_*g_*g^*\I, S^\dt(Y_\dt, f_*g_*g^*\I))\ar[r] & (f_*g_*g^*\I, h_*S^\dt(Y_\dt', {f'}_*g^*\I))
}
$$

\end{proof}    

\begin{lemma}\label{lemma:jouleray}
 Suppose that  $f:X\to Y$ is  a projective morphism in $V_\kappa$.
Let $h:Y'\to Y$ be a morphism with contractible fibers, 
and let $f':X'\to Y'$ denote the fiber product.
Then  the Leray spectral sequences
$$ H_\sigma^p(Y,R^qf_{\sigma*}\Z)\Rightarrow
H_\sigma^{p+q}(X,\Z)$$
and
$$ H_\sigma^p(Y',R^qf'_{\sigma*}\Z)\Rightarrow
H_\sigma^{p+q}(X',\Z)$$
are isomorphic.
\end{lemma}

\begin{proof}
There is of course, a map of spectral sequences induced by $h^*$.
It is enough to prove that this is an isomorphism at $E_2$.
We do this by induction on $dim\, Y$.
To simplify notation, we will omit $\sigma$ for the remainder
of the proof. 
Let $Y_1\subset Y$ be a closed set such that the sheaves $R^qf_*\Z$ 
are all locally constant on the complement of $Y_1$. Let $Y_1' =
h^{-1}Y_1$,
and let $i:Y_1\to Y$, $i':Y_1'\to Y'$, $j:Y-Y_1\to Y$ and
$j':Y'-Y_1'\to Y'$ denote the inclusions.
Then we have the following commutative diagram with exact rows:
$$
\xymatrix{
  H^{p-1}(Y_1,i_*R^qf_*\Z)\ar[r]\ar[d] & H^p(Y,j_!R^qf_*\Z)\ar[r]\ar[d] & H^p(Y,R^qf_*\Z)\ar[d] \\
 H^{p-1}(Y'_1,i'_*R^qf'_*\Z)\ar[r] &  H^p(Y',j'_!R^qf'_*\Z)\ar[r] &  H^p(Y',R^qf'_*\Z)
}
$$
The diagram should be extended two places to the right, but this is
omitted for typographic reasons.
The second and  fifth vertical arrows are isomorphisms because the pairs $(Y,Y_1)$
and $(Y',Y_1')$ with their classical topologies  are homotopy
equivalent. The first and fourth vertical arrows
are isomorphisms by induction.
Therefore the third vertical arrow is an isomorphism by the $5$ lemma.
\end{proof}

\begin{proof}[Proof of theorem \ref{thm:gen}]
Let  $f:X\to Y$ be  a projective morphism in $V_\kappa$.
By  lemma \ref{lemma:jouleray}, 
we can assume that $Y$ is affine.
Lemma \ref{lemma:existenceofcell}  yields a filtration $Y_\dt$ of $Y$
which is cellular with respect to the higher direct image sheaves  
$R^qf_{*}\Z$ (we will suppress $\sigma$) and satisfies $\dim Y_i= i$.
Let $X_\dt = f^{-1}Y_\dt$.
 We wish to show that the map of spectral sequences
 $$\lambda: H^p(Y, R^q f_*\Z) \to E_2^{pq}(X_\dt, \Z)$$
constructed in lemma ~\ref{lemma:leray} is an isomorphism.
It suffices to check that $\lambda$ is an  isomorphism at $E_2$.
Lemma \ref{lemma:cell} and  corollary \ref{cor:func2}
yield isomorphisms
$$H^p(Y,R^qf_*\Z)\cong E_2^{p0}(Y_\dt,R^q)
\cong E_2^{pq}(X_\dt, \Z) $$
We will be done once  we know that the composition of these
maps is $\lambda$. From the constructions of these maps, 
this amounts to the commutativity of the diagram
$$
\left.
\begin{array}{ccc}
 (\R f_* \Z, \tau) & \to & \R f_*(\Z,  Dec\, S(X_\dt)) \\
 \uparrow &  &\uparrow  \\
 (\tau_{\le q}\R f_* \Z, \tau) & \to  & (\tau_{\le q}\R f_*\Z, Dec\,  S(Y_\dt)) \\
 \downarrow &  & \downarrow \\
 (R^q[-q], \tau)  & \to & (R^q[-q], Dec\,  S(Y_\dt))
\end{array}
\right.
$$
This is a  straightforward verification.

It follows from  lemma~\ref{lemma:mhs} and the previous paragraph that
the Leray spectral sequence can be lifted to $\A$.
It remains to check that this lift is independent of $Y_\dt$
and functorial. Suppose that $Y_\dt'\subseteq Y$ is another chain  which is 
cellular with respect to the higher direct image sheaves and which satisfies 
the above dimension conditions.
Lemma \ref{lemma:existenceofcell} yields a third cellular filtration 
$Y_\dt''\supseteq Y_\dt\cup Y_\dt'$. Let $X_\dt' =
f^{-1}Y_\dt'$,  $X_\dt'' = f^{-1}Y_\dt''$ and let $\lambda',
\lambda''$ be the corresponding isomorphisms.
 Then we get  maps 
$$E_1^{pq}(X_\dt', \Z)\leftarrow E_1^{pq}(X_\dt'', \Z)\to E_1^{pq}(X_\dt, \Z)$$
compatible with $\lambda,\lambda'$ and $\lambda''$. Therefore, these are
isomorphisms from $E_2$ onwards. It follows that the liftings to $\A$ are isomorphic.

Next suppose that 
$$
\xymatrix{
 \tilde X\ar[r]^{g}\ar[d]^{\tilde f} & X\ar[d]^{f} \\
 \tilde Y\ar[r]^{h} & Y
}
$$
is a commutative diagram of varieties. By corollary \ref{cor:jou}
and lemma \ref{lemma:jouleray}, we can assume that both $Y$ and
$\tilde Y$ are affine. Let $\tilde Y_\dt\subseteq \tilde Y$ be a  filtration 
which is cellular with respect to all the  $R^q\tilde f_*\Z$ and such that
$\dim \tilde Y_i = i$. By lemma \ref{lemma:existenceofcell}, we can
construct a filtration $Y_\dt\supseteq h(\tilde Y_\dt)$ which
is cellular with respect to all $R^qf_*\Z$.
Let  $\tilde X_\dt = {\tilde f}^{-1}\tilde Y_\dt$. 
Then there is  a morphism of spectral sequences
$$E_2^{pq}(X_\dt, \Z)\to E_2^{pq}(\tilde X_\dt, \Z)$$
compatible with the map on Leray spectral sequences.
\end{proof}

\section{Compactly supported cohomology}

In this section, 
we will sketch a couple of generalizations of theorem~\ref{thm:gen}.
Once again, fix an embedding  $\sigma:\kappa\hookrightarrow \C$
and   an enriched $\sigma$-Betti theory
$(\cH:V_\kappa^2\to \A,\Phi)$.
Let $f:X\to Y$ be a projective morphism
of quasiprojective varieties over $\kappa$, and let $\tilde Y\subseteq Y$
be closed and $\tilde X = f^{-1}\tilde Y$. Let  $J:X-\tilde X\hookrightarrow  X$
and  $j:Y-\tilde Y\hookrightarrow  Y$ denote the inclusions.
There are isomorphisms $j_!R^qf_{\sigma,*}\Z =  R^qf_{\sigma,*}J_!\Z $.
Thus  we can form  the Leray  spectral sequence for pairs
\begin{equation}\label{eq:LerayPair}
\LL_2^{pq}(f,\tilde Y,\Z) = H_\sigma^p(Y, j_!R^qf_{\sigma*}\Z) \Rightarrow
H_\sigma^{p+q}(X,\tilde X,\Z)
\end{equation}
This is functorial in  the sense that a commutative diagram
$$
\xymatrix{
 (S,g^{-1}\tilde Z)\ar[r]^{g}\ar[d]^{f'} & (X,f^{-1}\tilde Y)\ar[d]^{f} \\
 (Z,\tilde Z)\ar[r]^{h} & (Y,\tilde Y)
}
$$
induces morphism 
$$
\LL_2^{pq}(f,\tilde Y,\Z)\to \LL_2^{pq}(f',\tilde Z, \Z).
$$

\begin{thm}
  This spectral sequence lifts naturally to $\A$.
\end{thm}

The proof is almost identical to the proof of theorem~\ref{thm:gen}, but a couple of
modifications  need to be made.   
We suppress $\sigma$ as above. We can reduce to the case where $Y$ is affine
by lemma~\ref{lemma:jou} and the following extension of
lemma~\ref{lemma:jouleray} to pairs:

\begin{lemma}
Let $h:Y'\to Y$ be a morphism with contractible fibers.
Let $f':X'\to Y'$ denote the fiber product
and $\tilde Y' = h^{-1}\tilde Y$ etcetera.
Then  the spectral sequences
$$ H^p(Y,j_!R^qf_{*}\Z)\Rightarrow
H^{p+q}(X,\tilde X,\Z)$$
and
$$ H^p(Y',j'_!R^qf'_{*}\Z)\Rightarrow
H^{p+q}(X',\tilde X',\Z)$$
are isomorphic.
\end{lemma}

This extension
follows immediately from lemma~\ref{lemma:jouleray} and the exact
sequence for a pair.
We then construct a cellular filtration $Y_\dt$ with respect to the 
sheaves $R^qf_{*}J_!\Z=j_!R^qf_{*}\Z$. Arguing as above, we see that 
 (\ref{eq:LerayPair}) is isomorphic to $E_2^{pq}(f^{-1}Y_\dt,J_!\Z)$.
The remaining ingredient is 
the following extension of lemma~\ref{lemma:mhs}:

\begin{lemma}
      For any pair $(X,X_\dt) $ in $FV_\kappa$ with
$J:X-\tilde X\to X$ as above, 
$$E_1^{pq}(X_\dt,J_!\Z)\Rightarrow H^{p+q}(X,J_!\Z)$$
can be lifted naturally to $\A$.
\end{lemma}

In fact, it is possible to use the original form of the lemma by 
choosing $\tilde Y$ as one of the members of the  cellular filtration;
this is possible since the sheaves $R^qf_{*}J_!\Z$
vanish along $\tilde Y$. The argument for independence from $Y_\dt$
and functoriality are essentially  as before.

\begin{cor}\label{cor:compactLeray}
Assume $\kappa=\C$. Let $f:X\to Y$ be a projective map of
quasiprojective varieties. Then the Leray spectral sequence for compactly supported
cohomology 
$$ \LL_{c,2}^{pq} = H_c^p(Y, R^qf_*\Z) \Rightarrow
H_c^{p+q}(X,\Z) 
$$
is compatible with mixed Hodge structures.
\end{cor}

\begin{proof}
  Choose projective compactifications $\bar X$ and $\bar Y$. By
  blowing up $\bar X$ along $\bar X- X$ we can assume that
$f$ extends to a morphism of the compactifications. Apply the
theorem to the map of pairs  
$(\bar X,\bar X- X)\to (\bar Y, \bar Y-Y)$.
\end{proof}

If $Y$ is affine we can give another proof of this corollary.

\begin{proof}[Second proof for $Y$ affine.]
  Proceed  as in the proof theorem~\ref{thm:gen}
and construct a cellular decomposition $Y_\dt$ with respect to
$R^qf_*\Z$, but now replace all occurrences of
ordinary cohomology with  compactly support cohomology.
In particular, the above spectral sequence will 
become isomorphic to
$$E_1 = H^{p+q}_c(X_p,X_{p+1}, \Z)\Rightarrow H^{p+q}_c(X,\Z)$$
from $E_2$ onwards, where $X_p= f^{-1}Y_p$.
\end{proof}

The second proof leads to an estimate on weights analogous
to the $l$-adic situation \cite{deligne3}.
Let us say that  a mixed Hodge structure $H$
has weights $\le k$ if $ Gr^W_iH=0$ for $i> k$.  

\begin{thm}\label{thm:weights}
  If $f:X\to Y$ is a projective map of complex quasiprojective
varieties, then $H_c^i(Y, R^jf_*\Z)$ has weights $\le i+j$.
\end{thm}

\begin{proof}
Suppose that $Y$ is affine, then from the previous proof
$H_c^i(Y, R^jf_*\Z)$  is a subquotient of
$ H^{i+j}_c(X_j,X_{j+1}, \Z)$, and this 
has weights $\le i+j$ by \cite[8.2.4]{deligne2}
and the exact sequence of the pair.

In general, choose an affine open set $U\subset Y$, and let
$Z= Y-U$. Then the theorem follows from induction on $\dim Y$
and the exact sequence
$$H_c^i(U, R^jf_*\Z)\to H_c^i(Y, R^jf_*\Z)\to H_c^i(Z, R^jf_*\Z)$$
\end{proof}

\begin{remark}
This proof yields something stronger. Namely, that the $(p,q)$, for which the
$(p,q)$th Hodge number of $H_c^i(Y, R^jf_*\Z)$ is nonzero,
lie within the shaded triangle of Figure 1.
\end{remark}

\begin{figure}[h]
\begin{center}
\includegraphics{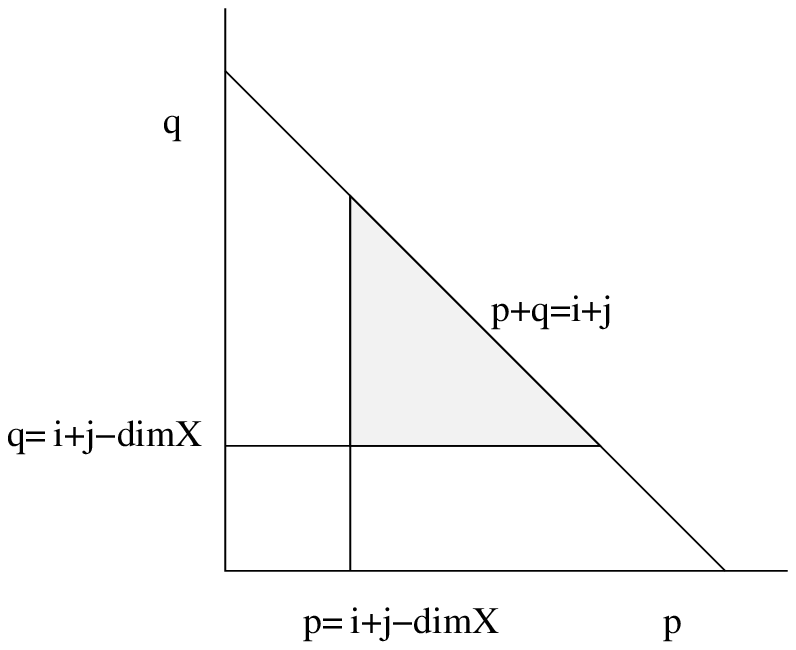}
\end{center}
\caption{}
\end{figure}

\section{Geometric Variations of Hodge Structure}

All varieties in this section will be defined over $\C$.
For the record, we will briefly recall the definition of
variation of Hodge structure; see \cite{bz}, and reference contained
there, for a more detailed discussion.
A rational polarizable variation of Hodge structure of weight $k$ 
on  a smooth  algebraic variety $Y$ consists of 
\begin{enumerate}
\item a local system or locally constant sheaf of $\Q$ vector spaces $V$ on $Y$.
\item a filtration by subbundles  $F^\dt \subset V\otimes O_X$ .
\end{enumerate}
We require that $F^\dt$ should satisfy  Griffiths transversality
and that the data induce a pure Hodge structure of weight $k$ on the fibers.
We also require  that the whole thing admits a flat polarization
(which will not be fixed). When no confusion is likely, we use the same
symbol for  the variation of Hodge structure and its local system.
A morphism is a homomorphism of local systems which preserves
filtrations.
Let $\PVHS_k(Y)$ be the additive (and in fact Abelian) 
category of polarizable variations
of Hodge structure over $Y$. The motivating examples come
from geometry. If $f:X\to Y$ is a smooth projective map of varieties
then $R^kf_*\Q$ is part of a polarizable variation of Hodge structure
of weight $i$ with
$$F^p = image[\R^k\Omega_{X/Y}^{\ge p}\to \R^k\Omega_{X/Y}^\dt\cong 
R^kf_*\Q\otimes O_Y]$$

 We  define a variation of Hodge structure to be {\em
geometric} if is isomorphic to a direct summand of some $R^if_*\Q$,
with $f$ smooth and projective. Let $\GVHS_k(Y)$ be the full subcategory
of geometric variations of Hodge structure. Given a morphism
$g:Y'\to Y$, there is an exact functor $g^*:\PVHS(Y)\to \PVHS(Y')$
defined by $V\mapsto g^*V$, $F^\dt\mapsto g^*F^\dt$.
$g^*$ preserves $\GVHS_k$, since $g^*R^kf_*\Q = R^kf'_*\Q$ where $f'$
is obtained by base change.

\begin{thm}\label{thm:VHS}
For each smooth quasiprojective variety $T$,
there are functors $H^i_{MHS}:\GVHS_k(T)\to \mbox{$\Q$-MHS}$ such that:
\begin{enumerate}
\item The rational lattice for $H^i_{MHS}(T,V)$ is just ordinary 
cohomology of the local system of $V$.
\item Let  $g:T'\to T$ be a morphism of varieties and let $V\in
  \GVHS(T)$, there is a morphism $g^*:H^i_{MHS}(T,V)\to H^i_{MHS}(T',g^*V)$
compatible with the usual pullback on cohomology.
\item If $f:X\to T$ is  a smooth projective map, then
$H_{MHS}^p(T, R^qf_*\Q)$ is a subquotient of the canonical mixed
Hodge structure on $H^{p+q}(X,\Q)$.
\end{enumerate}

\end{thm}

The proof will be based on a series of lemmas.

\begin{lemma}\label{lemma:multiplic}
If $F:Z\to T$ is a smooth projective morphism, then
$$L^pH^i(Z,\Q)\cup L^qH^j(Z,\Q)\subseteq L^{p+q}H^{i+j}(Z,\Q)$$
\end{lemma}

\begin{proof}
 Let $Z_p=F^{-1}T_p$ be as in theorem \ref{thm:smooth}. Then
the lemma follows from the commutativity of the diagram
$$
\xymatrix{
 H^i(Z)\otimes H^j(Z)\ar[r]\ar[dd] & H^i(Z_p)\otimes H^j(Z_q)\ar[d] \\
  & H^i(Z_{p+q})\otimes H^j(Z_{p+q})\ar[d] \\
 H^{i+j}(Z)\ar[r] & H^{i+j}(Z_{p+q})
}
$$

\end{proof}
Given a proper map $q:Z\to Y$ of smooth varieties.
The pushforward map $q_*:H^*(Z)\to H^*(Y)$ is Poincar\'e
dual to the pullback map $q^*$ on compactly
supported cohomology. This is a morphism of
mixed Hodge structure.

\begin{lemma}\label{lemma:pushforward}
 Suppose that
$$
\xymatrix{
 Z\ar[r]^{q}\ar[d]_{F} & Y\ar[ld]^{g} \\
 T & 
}
$$
is a commutative diagram of smooth projective maps
between smooth quasiprojective varieties.
Then
$$q_*(L^pH^i(Z,\Q))\subseteq L^pH^{i-2d}(Y,\Q)$$
where $d= \dim Z - \dim Y$.
\end{lemma}

\begin{proof}
Let  $Y_p = g^{-1}T_p$ and $Z_p = F^{-1}T_p= q^{-1}Y_p$ in the notation of
 theorem \ref{thm:smooth}.  We have a commutative diagram
$$
\xymatrix{
 H^i(Z)\ar[r]^{q_{*}}\ar[d] & H^{i-2d}(Y)\ar[d] \\
 H^i(Z_p)\ar[r]^{q_{*}} & H^{i-2d}(Y_p)
}
$$
The lemma follows from this and   theorem \ref{thm:smooth}.
\end{proof}

Let $f:X\to T$ and $g:Y\to T$ be smooth projective
maps. Let us write $f\times g:X\times_T Y \to T$ for the
fiber product, and $p_1:X\times_T Y \to X$ and
$p_2:X\times_T Y \to Y$ for the projections.
 Suppose that the fiber dimension of $f$ is $d$.
Given a morphism $\gamma:R^if_*\Q \to R^ig_*\Q $,
the class of its graph $[\gamma]$ is the image
of $\gamma$ under
\begin{eqnarray*}
 Hom_{VHS}(R^if_*\Q, R^ig_*\Q) &\cong&
Hom_{VHS}(\Q(0), R^{2d-i}f_*\Q(d)\otimes R^ig_*\Q)\\
&\subseteq& Hom_{VHS}(\Q(0), R^{2d}(f\times g)_*\Q(d))
\end{eqnarray*}
Thus $[\gamma]$ can be viewed as a constant degree $2d$ 
Hodge cycle, i.e. rational $(d,d)$ cycle, along the fibers of $f\times g$.

\begin{lemma}\label{lemma:corresp}
With the above notation, $\gamma$ coincides with  the composite of
$$R^if_*\Q \stackrel{p_1^*}{\longrightarrow}
R^{i}(f\times g)_*\Q \stackrel{\cup [\gamma]}{\longrightarrow}
R^{i+2d}(f\times g)_*\Q(d) \stackrel{p_{2*}}{\longrightarrow}
R^if_*\Q
$$
\end{lemma}    

\begin{proof}
It is enough to check this at the stalk level, and this is
well known \cite[1.3.7]{kleiman}.  
\end{proof}

\begin{proof}[Proof of theorem \ref{thm:VHS}]
By \cite{deligne}, there is an isomorphism of vector spaces
$$H^p(T, R^qf_*\Q)\cong Gr_L^pH^{p+q}(X,\Q)$$
Corollary \ref{cor:MHSleray} implies that the right hand side
carries a mixed Hodge structure. We use this isomorphism 
to define the mixed Hodge structure on the left. Thus (3) is
true by decree. More generally, suppose that $V$  is a variation
of Hodge structure which is a direct summand of $R^if_*\Q$.
Then $V$ is the image of a projection
$P:R^if_*\Q\to  R^if\Q $.
By the theorem of the fixed part \cite[4.1.1]{deligne2},
the graph 
$$[P]\in Hom_{VHS}(\Q(0), R^{2d}(f\times f)_*\Q(d))$$
can be lifted to a Hodge cycle $\phi\in H^{2d}(X,\Q)$. 
The map $\Phi(\alpha)=p_{2*}(p_1^*\alpha\cup \phi)$ induces a
morphism
$\Phi:H^{i+j}(X,\Q)\to H^{i+j}(X,\Q)$
of mixed Hodge structures, which by 
(\ref{eq:basechangeLeray}) and lemmas \ref{lemma:multiplic} and
 \ref{lemma:pushforward}, 
 preserves the Leray filtrations. By  lemma \ref{lemma:corresp},
$$Gr_L^i(\Phi):H^j(X,R^if_*\Q)\to H^j(X,R^if_*\Q)$$
coincides with the map induced by $P$. This shows that
$P$ is a morphism of mixed Hodge structures. 
Therefore we can equip $H^j(X,V)$ with the mixed Hodge structure on 
the image  $P(H^j(T,R^if_*\Q))$.

We have to verify that the mixed Hodge structure just constructed on
$H^j(X,V)$ is independent of choices.  Suppose that $V$ is a direct summand
of both $R^if_*\Q$ and $R^ig_*\Q$ for two smooth projective maps
$f:X\to T$ and $g:Y\to T$. Then consider the commutative
diagram
$$
\xymatrix{
 R^if_*\Q\ar[rr]^{P}\ar[rd]\ar@/^/[rrdd]^{\gamma} &  & R^if_*\Q \\
  & V\ar[rd]\ar[ru] &  \\
 R^ig_*\Q\ar[rr]_{\Pi}\ar[ru]\ar@/_/[rruu]_{\delta} &  & R^ig_*\Q
}
$$
where all the unlabeled arrows  are the projections and
inclusions to and from $V$. $\gamma$ ($\delta$) is
the arrow connecting the upper (respectively lower) left and lower
(respectively upper) right corners.
The relations $\gamma P = \gamma$ and
$\gamma \delta =\Pi$ are easily verified.
We get an induced  diagram
$$
\xymatrix{
 H^i(X,R^jf_*\Q)\ar[r]^{P}\ar@/^/[rd]^{\gamma} & H^i(X,R^jf_*\Q) \\
 H^i(Y,R^jg_*\Q)\ar[r]_{\Pi}\ar@/_/[ru]_{\delta} & H^i(Y,R^jg_*\Q)
}
$$
We can see that these maps are morphisms of mixed Hodge structures
by expressing them as   compositions of
pullbacks, pushforwards and cup products with Hodge cycles
as above.
The relations $\gamma P = \gamma$ and $\gamma \delta =\Pi$
imply that $\gamma$ induces a surjection
$image(P)\to image(\gamma)$ and that there is an inclusion 
$image(\Pi)\subseteq image(\gamma)$.
Since 
$$\dim image(P) = \dim image(\Pi)= \dim H^i(X,V),$$
it follows that $image(P)$ and $image(\Pi)$ are isomorphic
as mixed Hodge structures.

Finally property (2) is clear from construction.

\end{proof}

\begin{cor}
Let $f:X\to T$ be a smooth projective map of smooth
quasiprojective varieties, 
then $H^j(T,S^N(R^if_*\Q))$ carries a natural mixed Hodge
structure for each positive integer $N$.
The same holds if the symmetric power is replaced by 
an exterior power or,  more generally, a  Schur functor. 
\end{cor}

\begin{proof}
  These are  geometric variations of Hodge structure, 
since they  are direct summands of 
$R^i(f\times f\times \ldots f)_*\Q$.
\end{proof}

\section{Threefolds fibered by rational curves}

We end this paper with a geometric example. The details
are sketched.
Let $X$ be a smooth projective threefold over $\C$ with
a flat map $f:X\to S$ such that the fibers are connected
rational curves.  Let $C$ be the discriminant set. 
Our goal is to compute the intermediate
Jacobian 
$$J^2(X) = \frac{H^3(X,\C)}{F^2H^3(X,\C)+H^3(X,\Z)}$$
generalizing the usual description for conic bundles \cite{cg}.
This amounts to computing the Hodge structure on $H^3(X)$
modulo torsion. 

We start by computing higher direct images.
The fibers of $f$ are rational curves which are homotopic
to a wedge of spheres. Therefore $f_*\Z =\Z$, 
$$(R^2f_*\Z)_s\cong \Z^{(\mbox{$\#$ irreducible components of
    $f^{-1}(s)$})}$$
and all other $R^if_*\Z=0$. Thus the Leray spectral sequence
(\ref{cor:LerayMHS}) yields an exact sequence of Hodge structures
\begin{equation}\label{eq:H3seq}
0\to \LL^{30}_\infty\to H^3(X,\Z) \to  \LL^{12}_\infty\to 0
\end{equation}
We certainly have $\LL^{30}_2=\LL^{30}_3$ and $\LL^{12}_2=\LL^{12}_3$.
Consider the differentials
$$d_3:\LL^{02}_3\to \LL^{30}_3$$
$$d_3:\LL^{12}_3\to \LL^{40}_3$$
The right sides are pure of weights $3$ and $4$ respectively, while the
left sides have smaller weights by theorem \ref{thm:weights}.
Therefore these differentials vanish. (This can also be deduced
by observing that
the maps $H^i(S,\Z)\to H^i(X,\Z)$ are injective modulo
torsion because of Poincar\'e duality.) Thus
$\LL^{30}_\infty = H^3(S,\Z)$ and  $\LL^{12}_\infty = H^1(S,R^2f_*\Z)$.
The sequence (\ref{eq:H3seq}) leads an exact sequence
of tori
$$0\to Alb(S)\to J^2(X)\to J^2 H^1(S,R^2f_*\Z)\to 0$$
To proceed further, let us make the simplifying assumptions
that $Alb(S)=0$ or equivalently that  $H^1(S,\Z)$ is torsion,
that $C$ is smooth, and that $f^{-1}C \to C$ is topologically
locally trivial.
Choose an embedding 
$$
\xymatrix{
 X\ar[d]^{f}\ar[r]^{i} & \PP^N\times S\ar[ld]^{g} \\
 S & 
}
$$
Then $R^2g_*\Z\cong\Z_S$ injects into $R^2f_*\Z$. The
 torsion free part $L$ of $R^2f_*\Z/R^2g_*\Z$ is supported
on $C$. From the long exact sequence, we obtain an isomorphism
$$H^1(S,R^2f_*\Z)\cong H^1(C,L)$$
Consider the diagram
$$
\xymatrix{
   Z\ar[r]\ar[d] & f^{-1}C\ar[d] \\
   \tilde C\ar[r]^{\pi} & C
}
$$
where the $Z\to \tilde C$  is the Stein factorization 
of normalization of $f^{-1}C$.
$L\otimes \Q$ is a quotient of $\pi_*\Q$.
It follows easily that $J^2(X)$ is isogenous to a quotient
of the Jacobian $J(\tilde C)$. 
In the case of conic bundle,
 $L$ is the anti-invariant part of  $\pi_*\Z$
under the involution of the double cover $\tilde C\to C$.
So we recover the usual description of $J^2(X)$ as a Prym
variety.


\begin{thebibliography}{AB}

\bibitem[AS]{as} D. Arapura, P. Sastry, {\em Hodge structures
and intermediate Jacobians of moduli spaces}, Proc. Ind.
Acad. Sci 110 (2000)

\bibitem[BBD]{bbd} A. Beilinson, J. Bernstein, P. Deligne,
{\em Faiseaux Perverse}, Asterisque 100 (1982)



\bibitem[BZ]{bz} J. L. Brylinski, S. Zucker, {\em An
overview of recent advances in Hodge theory}, Sev. Complex
Variables VI, Springer-Verlag (1990)

\bibitem[CG]{cg}C. H. Clemens, P. Griffiths
{\em The Intermediate Jacobian of the Cubic Threefold} 
Ann.  Math. 95 (1972)

\bibitem[D1]{deligne} P. Deligne, {\em Th\'eor\`eme de Lefschetz
et crit\'eres de d\'eg\'en\'erescence de suites spectrales}
Publ. IHES 35 (1968)

\bibitem[D2]{deligne2} P. Deligne, {\em Th\'eorie de Hodge II, III}
Publ. IHES 40, 44, (1971, 1974)

\bibitem[D3]{deligne3} P. Deligne, {\em La conjecture de Weil II},
Publ. IHES 52 (1980)

\bibitem[H]{hain} R. Hain,
{\em Torelli groups and geometry of moduli spaces of curves.}
Current topics in complex geometry, MSRI Publ. (1995)


\bibitem[HL]{hamm} H. Hamm,  D.T. Le,
{\em Vanishing theorems for constructible sheaves. II.}
 Kodai Math. J.  21  (1998)

\bibitem[GM]{gm} M. Goresky, R. Macpherson, {\em Intersection
homology II}, Invent. Math 72 (1983)

\bibitem[J]{jannsen} U. Jannsen, {\em Mixed motives and
algebraic K-theory}, Lect notes in math 1400, Springer-Verlag (1990)

\bibitem[Jo]{jou} J.P. Jouanolou, {\em Un suit exacte de
    Mayer-Vietoris en K-theorie algebrique},in Algebraic K-theory,
Lect. Notes Math 341, Springer-Verlag (1973)


\bibitem[K]{kleiman} S. Kleiman, {\em Algebraic cycles and the
Weil conjectures}, Dix Esposes sur la cohomologies les schemas,
North Holland (1968)

\bibitem[L]{lewis} J. Lewis, {\em A filtration on the
Chow groups of a complex projective variety}, Compositio Math. 128
(2001)



\bibitem[N]{nori} M. Nori, {\em Constructible sheaves}, Proc. Int.
Conf. on Algebra..., TIFR (2002)


\bibitem[P]{peters} C. Peters, {\em Letter to the author}, May 26, 2004

\bibitem[S1]{saito1} M. Saito, {\em Modules Hodge polarizables}
Publ. RIMS 24 (1988)

\bibitem[S2]{saito2} M. Saito, {\em Mixed Hodge modules}
Publ. RIMS 26 (1990)

\bibitem[W]{weibel} C. Weibel, {\em An introduction to
homological Algebra}, Cambridge U. Press (1994)


\bibitem[Z]{zucker} S. Zucker, {\em Hodge theory with 
degenerating coefficients...}, Annals of Math. 109 (1979)


\end{thebibliography}
\end{document}